\documentclass[aap,MSNbibl,seceqn,citesort,dvips]{arximspdf}

\doi{10.1214/09-AAP638}
\volume{20}
\issue{3}
\pubyear{2010}
\firstpage{869}
\lastpage{889}

\makeatletter
\newtheorem{theorem}{Theorem}
\newtheorem{lemma}[theorem]{Lemma}
\newtheorem{prop}[theorem]{Proposition}

\newproclaim{remark}{Remark}
\newproclaim{notation}{Notation}

\def\fract{\tfrac}
\def\eqref#1{(\ref{#1})}
\makeatother

\begin{document}
\begin{frontmatter}

\title{The random conductance model with Cauchy tails\protect\thanksref{T1}}
\runtitle{RCM with Cauchy tails}

\begin{aug}
\author[A]{\fnms{Martin T.} \snm{Barlow}\ead[label=e1]{barlow@math.ubc.ca}} \and
\author[B]{\fnms{Xinghua} \snm{Zheng}\corref{}\ead[label=e2]{xhzheng@ust.hk}}
\runauthor{M. T. Barlow and X. Zheng}

\thankstext{T1}{Supported in part by NSERC (Canada) and the Peter Wall Institute for Advanced Studies.}
\affiliation{University of British Columbia and\break Hong Kong
University of Science and Technology}
\address[A]{Department of Mathematics\\
University of British Columbia\\
Vancouver, British Columbia V6T 1Z2\\
Canada\\
\printead{e1}} 
\address[B]{Department of ISOM\\
Hong Kong University of Science\\
\quad and Technology\\
Clear Water Bay, Kowloon\\
Hong Kong\\
\printead{e2}}
\end{aug}

\received{\smonth{8} \syear{2009}}

%
\begin{abstract}
We consider a random walk in an i.i.d. Cauchy-tailed
conductances environment. We obtain a quenched functional CLT for
the suitably rescaled random walk, and, as a key step in the
arguments, we improve the local limit theorem for $p^{\omega}_{n^2
t}(0,y)$ in [\textit{Ann. Probab.} (2009). To appear], Theorem~5.14,
to a result which gives
uniform convergence for $p^{\omega}_{n^2 t}(x,y)$ for all $x, y$ in a
ball.
\end{abstract}

%


\end{frontmatter}

\setcounter{section}{-1}
\section{Introduction}

In this paper we will establish the convergence to Brownian motion
of a random walk in a symmetric random environment in a critical
case that has not been covered by the papers \cite{BC09,BD08}. We
begin by recalling the ``random conductance model'' (RCM). We
consider the Euclidean lattice $\mathbb{Z}^d$ with $d\ge2$. Let $E_d$ be
the set of nonoriented nearest neighbour bonds, and, writing
$e=\{x,y\} \in E_d$, let $(\mu_e, e \in E_d)$ be nonnegative
i.i.d. r.v. on $[1,\infty)$ defined on a probability space
$(\Omega, \mathbb{P})$. We write $\mu_{xy}=\mu_{\{x,y\}}=\mu
_{yx}$; let
$\mu_{xy}=0$ if $x \not\sim y$, and set $\mu_x =\sum_y \mu_{xy}$.

We consider two continuous time random walks on $\mathbb{Z}^d$ which jump
from $x$ to $y\sim x $ with probability $\mu_{xy}/\mu_x$. These
are called in \cite{BD08} the \textit{constant speed random walk}
(CSRW) and \textit{variable speed random walk} (VSRW), and have
generators
%
\begin{eqnarray}\label{e-LC}
\mathcal{L}_C({\omega}) f(x) &=& \mu_x({\omega})^{-1} \sum_y \mu_{xy}({\omega})\bigl(f(y)-f(x)\bigr),\\\label{e-LV}
\mathcal{L}_V({\omega}) f(x) &=& \sum_y \mu_{xy}({\omega}) \bigl(f(y)-f(x)\bigr).
\end{eqnarray}
We write $X$ for the CSRW and $Y$ for the VSRW. Thus $X$ jumps
out of a state $x$ at rate~1 while $Y$ jumps out at rate $\mu_x$.
We will abuse notation slightly by writing $P^x_{\omega}$ for the laws
of both $X$ and $Y$ started at $x \in\mathbb{Z}^d$ in the random
environment $[\mu_e({\omega})]$. Since the generators of these
processes differ by a multiple, $X$ and $Y$ are time changes of
each other. More explicitly, as in \cite{BC09}, define the \textit{clock process}
%
\begin{equation}\label{e:defclock}
S_t = \int_0^t \mu_{Y_s}\,ds,
\end{equation}
and let $A_t$ be its inverse. Then the CSRW can be defined by
%
\begin{equation}\label{e:defcsrw}
X_t = Y_{A_t},  \qquad t\ge0.
\end{equation}

In the case when $\mu_e \in[0, 1]$, and $\mathbb{P}(\mu_e>0)> p_c(d)$,
the critical probability for bond percolation in $\mathbb{Z}^d$, the
papers \cite{BP07,Mat08} prove that both $X$ and $Y$ satisfy a
quenched functional central limit theorem (QFCLT), and that the
limiting process is nondegenerate. The paper \cite{BD08} studies
the case when $\mu_e \in[1, \infty)$, and proves that for
$\mathbb{P}$-a.a. ${\omega}$ the rescaled VSRW, defined by
%
\begin{equation}\label{e:defyn}
Y^{(n)}_t = n^{-1} Y_{n^2t}, \qquad  t \ge0,
\end{equation}
converges to $(\sigma_V W_t, t \ge0)$ where $W$ is a standard
Brownian motion, and $\sigma_V>0$. It is also proved there that
$S_t/t \rightarrow \mathbb{E}\mu_0\in[1, \infty]$. It follows from
\eqref{e:defcsrw} that the CSRW with the standard rescaling,
\[
X^{(n,1)}_t = n^{-1} X_{n^2t}, \qquad  t \ge0,
\]
converges to $\sigma_C W$ where
\begin{eqnarray*}
\sigma_C =
\cases{
\sigma_V / \sqrt{2d  \mathbb{E}\mu_e}, &\quad \mbox{if } $\mathbb{E}\mu_e<\infty$,\cr
0, &\quad \mbox{if } $\mathbb{E}\mu_e = \infty$.
}
\end{eqnarray*}
If $\mathbb{E}\mu_e =\infty$ it is natural to ask if a different
rescaling of $X$ will give a nontrivial limit. In the case when
$d\ge3$, $\mu_e \in[1,\infty)$ and there exists $\alpha\in
(0,1)$ such that
%
\begin{equation}\label{e:mutail}
\mathbb{P}(\mu_e > u ) \sim\frac{c}{ u^{\alpha}}\qquad  \mbox{as }u\rightarrow\infty,
\end{equation}
then \cite{BC09} proves that the process
\[
X^{(n,\alpha)}_t = n^{-1} X_{n^{2/\alpha}t}, \qquad  t \ge0,
\]
converges to the ``fractional kinetic motion'' with index $\alpha$.
(For details of this process, and its connection with aging see
\cite{BCM06,BC07,BC08}.) These papers leave open the case when
$\alpha=1$. In this paper we assume that $(\mu_e)$ satisfies
\eqref{e:mutail} with $\alpha=1$; for simplicity we take
$c=1/(2d)$, so that $\mu_e$ satisfies
%
\begin{eqnarray}\label{e:muh1}
\mathbb{P}(\mu_e \ge1) &=& 1, \\\label{e:muxtail}
\mathbb{P}(\mu_e\geq u ) &\sim&\frac{1}{2du}\qquad  \mbox{as }
u\rightarrow \infty.
\end{eqnarray}

We define the process
%
\begin{equation}\label{e:defXn}
X^{(n)}_t = n^{-1} X_{n^2 (\log n)  t},\qquad   t \ge0.
\end{equation}

Our main theorem follows:
\begin{theorem}\label{thm:main} Let $d \ge3$, and assume that $\mu_e$ satisfies
\eqref{e:muh1} and \eqref{e:muxtail}. Then for $\mathbb{P}$-a.a.
${\omega}$, $(X^{(n)}, P^0_{\omega})$ converges in $D([0,\infty
);\mathbb{R}^d)$
to $\sigma_1 W$ where $\sigma_1= \sigma_V/\sqrt{2} >0$, and $W$
is a standard $d$-dimensional Brownian motion.
\end{theorem}

As in \cite{BC09} we prove this theorem by using \eqref{e:defcsrw}
and proving convergence of a rescaled clock process. Let
%
\begin{equation}
S^{(n)}_t = \frac{1}{n^2\log n} \int_0^{n^2 t} \mu_{Y_s}\,ds;
\end{equation}
then it is easy to check that if $A^{(n)}$ is the inverse of $S^{(n)}$,
then
%
\begin{equation}
X^{(n)}_t = Y^{(n)}_{A^{(n)}_t}, \qquad  t\ge0.
\end{equation}
It follows that to prove Theorem~\ref{thm:main} it is enough to
prove.
\begin{theorem}\label{thm:conv}
Let $d \ge3$, and assume that $\mu_e$ satisfies \eqref{e:muh1}
and \eqref{e:muxtail}. For $\mathbb{P}$-a.a. ${\omega}$,
under the law
$P^0_{\omega}$,
%
\begin{equation}\label{eq:conv_clock}
\bigl(S^{(n)}_t, t \ge0\bigr) \Rightarrow(2t, t\ge0)\qquad \mbox{on } C([0,\infty);\mathbb{R}).
\end{equation}
\end{theorem}
\begin{remark}
For $\lambda\in[1, \infty)$, let
$S_t^{(\lambda)}=\frac{1}{\lambda^2\log\lambda}
\int_0^{\lambda^2t} \mu_{Y_s}\,ds$. Then if $n\leq\lambda\leq
(n+1)$,
\[
\frac{n^2\log n}{(n+1)^2\log(n+1)}\cdot S_t^{(n)} \le S_t^{(\lambda
)} \le
\frac{(n+1)^2\log(n+1)}{n^2\log n }\cdot S_t^{(n+1)}.
\]
It follows that the convergence \eqref{eq:conv_clock} holds for
$(S_t^{(\lambda)}, t\geq0)_{\lambda\ge1}$, and hence Theorem~\ref{thm:main}
extends to $(X_t^{(\lambda)})_{\lambda\ge 1}:=(\lambda^{-1} X_{\lambda^2 (\log\lambda) t})_{\lambda\ge 1}$.
\end{remark}

As in \cite{BC09}, the result is proved by estimating the growth
of the clock process $S_t$, $0\le t \le n^2 T$. Since the limit of
the processes $S^{(n)}$ is deterministic, overall this case is much
easier than when $\alpha\in(0,1)$: after suitable truncation it
is enough to use a mean--variance calculation. There is, however, one
respect in which this case is more delicate than when $\alpha<1$.
When $\alpha<1$ it turns out that the main contribution to
$S_{n^2 T}$ is from visits by $Y$ to $x$ such that $\varepsilon
n^{2/\alpha }
\le\mu_x \le\varepsilon ^{-1} n^{2/\alpha }$ (see Sections 5 and 7 of
\cite{BC09}). When $\alpha=1$ one finds that each set of edges of
the form $E_i=\{e{}\dvtx{} 2^{i-1} n \le\mu_e < 2^{i} n\}$, $i=1, \ldots,
\log n$, has a roughly comparable contribution to $S_{n^2 T}$, so
a much greater range of values of $\mu_e$ need to be considered.

To motivate the proof, consider the classical case of a sum of
i.i.d. r.v. $\xi_i$, with $\mathbb{P}( \xi_i > t) \sim
t^{-1}$. We
have that if
%
\begin{equation}\label{e:class1}
U^{(n)}_t = (n \log n)^{-1} \sum_{i=1}^{[nt]} \xi_i,
\end{equation}
then\vspace*{1pt} $\sup_{0\le t\le T} |U^{(n)}_t - t| \rightarrow 0$ in probability.
Let $a_i = i (\log i)^\beta$ where $\beta\in(1,2)$, and $\xi'_i =
\xi_i \mathbf{1}_{(\xi_i
> a_i)}$. Then $\sum P( \xi_i \neq\xi'_i)$ converges, so it is
enough to consider the convergence of
%
\begin{equation}\label{e:class2}
V^{(n)}_t = (n \log n)^{-1} \sum_{i=1}^{[nt]} \xi'_i.
\end{equation}
A straightforward argument calculating the mean and variance of
%
\begin{equation}\label{e:class3}
M^{(n)}_t = (n \log n)^{-1} \sum_{i=1}^{[nt]}( \xi'_i - E \xi_i')
\end{equation}
then gives convergence of $U^{(n)}$.
[Note that one does not have a.s. convergence, since
$P(\max_{2^{n-1} \le i \le2^n} \xi_i > 2^n \log2^n) \sim c/n$.]

The equivalent arguments in our case rely on good control of the
process~$Y$. Define the heat kernel and Green's functions for $Y$
by
%
\begin{equation}
p^{\omega}_t(x,y) = P^x_{\omega}(Y_t=y),\qquad
g^{\omega}(x,y) =\int_0^\infty p^{\omega}_t(x,y)\, dt.
\end{equation}
We extend these functions from $\mathbb{Z}^d \times\mathbb{Z}^d$ to
$\mathbb{R}^d \times
\mathbb{R}^d$
by linear interpolation on each cube in $\mathbb{R}^d$ with vertices
in $\mathbb{Z}^d$.
Let $W$ be a standard Brownian motion on $\mathbb{R}^d$, and let
$W^*_t =
\sigma_V W_t$,
so that $W^*$ is the weak limit of the processes~$Y^{(n)}$. Let
%
\begin{equation}\label{e:def_k}
k_t(x) = (2\pi\sigma^2_V)^{-d/2} \exp( -|x|^2/2 \sigma_V^2)
\end{equation}
be the density of the $W^*$.

A key element of the arguments is the following strengthening of
the local limit theorem for $p^{\omega}_{n^2 t}(0,y)$ in \cite{BD08}, Theorem~5.14,\vspace*{1.5pt} to a result which gives uniform convergence for
$p^{\omega}_{n^2 t}(x,y)$ for all $x, y$ in a ball.
\begin{theorem}\label{thm:unif_llt}
Let $d \ge2$, and assume $\mu_e$ satisfies \eqref{e:muh1}. For
any $\varepsilon>0$, $0<\delta<T<\infty$ and $K> 0$, we have the
following $\mathbb{P}$-almost sure \textit{uniform} convergence:
\begin{eqnarray}\label{eq:tp_unif_conv}
\nonumber
\frac{1}{1+\varepsilon} &<&
\liminf_{n\rightarrow \infty}\inf_{\delta\leq t\leq T}\inf_{|x|,|y| \leq K}
\frac{n^d p_{n^2t}^\omega(nx,ny)}{k_{t}(x,y)} \\[-8pt]\\[-8pt]
&\le&\limsup_{n\rightarrow \infty}\sup_{\delta\leq t\leq T}\sup_{|x|,|y|
\leq K} \frac{n^d p_{n^2t}^\omega(nx,ny)}{k_{t}(x,y)} < 1+\varepsilon.\nonumber
\end{eqnarray}
\end{theorem}

This result is proved in Section~\ref{ssec:unif_llt}.
\begin{notation*}
We write
\[
B(x,r) =\{ y \in\mathbb{Z}^d{}\dvtx{} |x-y| \le r\}\quad  \mbox{and}\quad
B_\mathbb{R}(x,r) =\{ y \in\mathbb{R}^d{}\dvtx{} |x-y| \le r\}.
\]
If $e=\{x_e, y_e\} \in E_d$, we write $e \in B(x,r)$ if $\{x_e,
y_e\} \subset B(x,r)$. We will follow the custom of writing $f\sim
g$ to mean that the ratio $f/g$ converges to 1, and $f\asymp g$ to
mean that the ratio $f/g$ remains bounded away from 0
and~$\infty$. For any $a,b\in\mathbb{R}$, $a\wedge b:=\min(a,b)$, and
$a\vee b:=\max(a,b)$. Throughout the paper, $c, C,C_1, C'$,
et cetera, denote generic constants whose values may change from
line to line.
\end{notation*}
\begin{remark}
One can also consider the more general case when the tail of
$\mu_e$ satisfies
\[
\mathbb{P}(\mu_e\geq u )\sim c\frac{(\log u)^\rho}{u}\qquad \mbox{as }
u\rightarrow\infty,
\]
where $\rho\geq-1$ (so that $\mathbb{E}\mu_e=\infty$). Define for
$t \ge
0$
\begin{eqnarray*}
X^{(n)}_t =  \cases{ n^{-1} X_{n^2 (\log n)^{1+\rho} t}, &\quad\mbox{when } $\rho> -1$,\cr
n^{-1} X_{n^2 (\log\log n)  t},&\quad\mbox{when } $\rho= -1$.}
\end{eqnarray*}
Then using the same strategy as in this article one can show that
for $\mathbb{P}$-a.a.~${\omega}$, $(X^{(n)},P^0_{\omega})$ converges
to a (multiple of a) Brownian motion.
\end{remark}

\section{Preliminaries}\label{sec:prelim}

\subsection[Heat kernel: Proof of Theorem 3]{Heat kernel: Proof of Theorem~\protect\ref{thm:unif_llt}}\label{ssec:unif_llt}

We collect some known estimates for $p^{\omega}_t(x,y)$ and
$g^{\omega}(x,y)$ which will be used in our arguments.
\begin{lemma}\label{lem:hk_est}
Let $\eta\in(0,1)$. There exist random variables $U_x$ $(x\in
\mathbb{Z}^d)$ and constants $c_i$ such that
\begin{eqnarray*}
\mathbb{P}(U_x\geq n)\leq c_1\exp(-c_2n^\eta),\qquad \mbox{for all } n\ge 1.
\end{eqnarray*}

\begin{longlist}[(a)]
\item[(a)]\label{item:tp_unif_ub} \cite{BD08}, Theorem~1.2(\textup{a}). There
exists $c_3>0$ such that for all $x,y$ and $t$,
\begin{eqnarray*}
p^\omega_{t}(x,y) \leq c_3 t^{-d/2}.
\end{eqnarray*}
\item[(b)]\label{item:tp_ub} \cite{BD08}, Theorem~1.2(\textup{b}). If
$|x-y|\vee\sqrt{t}\geq U_x$, then
%
\begin{eqnarray}\label{eq:hk_bd}
\hspace*{30pt}&&p_t^\omega(x,y)\nonumber
\\[-8pt]\\[-8pt]
&&\qquad\leq
\cases{
 c_4 t^{-d/2}\exp(-c_5 |x-y|^2/t), &\quad \mbox{when } $t\geq|x-y|$,\cr
 c_4\exp \bigl(-c_5 |x-y|\bigl(1 \vee\log(|x-y|/t) \bigr)  \bigr),&\quad\mbox{when } $t\leq|x-y|$.
 }\nonumber
\end{eqnarray}
\item[(c)]\label{item:tp_lb} \cite{BD08}, Theorem~1.2(\textup{c}). If $t\geq
U_x^2\vee|x-y|^{1+\eta}$, then
\[
p_t^\omega(x,y) \geq c_6 t^{-d/2}\exp(-c_7|x-y|^2/t).
\]
\item[(d)] Let $\tau(x,R) = \inf\{t\ge0{}\dvtx{} |Y_t-x|>R\}$. If $R \ge U_x$,
then
\[
P^x_{\omega}\bigl( \tau(x,R) \le t\bigr) \le c_{8}\exp( -c_{9} R^2/t ).
\]
\item[(e)]\cite{BC09}, Lemma~3.4. When $d\geq3$,
%
\begin{equation}\label{e:Gn_bnd}
c_{10} U_x^{2-d}\leq g^\omega(x,x)\leq c_{11}.
\end{equation}
\item[(f)]\cite{BC09}, Proposition~3.2(\textup{b}). When $d\geq3$, if $|x|\geq
U_0$, then
%
\begin{equation}\label{eq:Gn_gen_bd}
g^\omega(0,x)\leq\frac{c_{12}}{|x|^{d-2}}.
\end{equation}
\item[(g)] \cite{BC09}, Lemma~3.3.
There exists $c_{13}>0$ such that for each $K>0$, if
%
\begin{equation}\label{e:b1-def}
b_n=c_{13} (\log n)^{1/\eta},
\end{equation}
then with $\mathbb{P}$-probability no less than $1-c_{14} K^d n^{-2}$ the
following holds:
%
\begin{equation}\label{eq:U_x}
\max_{|x|\leq Kn} U_x\leq b_n.
\end{equation}
In particular, \eqref{eq:U_x} holds for all $n$ large enough
$\mathbb{P}$-a.s.
\item[(h)] \cite{BD08}, Theorem~5.14. For any
$\delta>0$, $\mathbb{P}$-a.s.,
%
\begin{equation}\label{e:llt}
\lim_{n\rightarrow \infty} \sup_{x\in\mathbb{Z}^d} \sup_{t\geq
\delta}
 | n^d p^\omega_{n^2 t} (0,x) - k_t (x/n) | = 0.
\end{equation}
\item[(i)] There exists $\theta>0$ such that for $x,
y, y' \in\mathbb{Z}^d$,
%
\begin{equation}\label{e:holder}
n^d  |p^\omega_{n^2t}(x,y) - p^\omega_{n^2t}(x,y') | \leq
c_{15} t^{-(d+\theta)/2} \cdot \biggl( \frac{ |y-y'| \vee U_y} {n} \biggr)^\theta.
\end{equation}
\end{longlist}
\end{lemma}

\begin{pf}
(d) The tail bound on $\tau(x,R)$ in (d) follows from Proposition~2.18 and Theorem~4.3 of \cite{BD08}.
(i) This follows from \cite{BD08}, Theorem~3.7 and \cite{BH09}, Proposition~3.2.
\end{pf}

We begin by improving the local limit theorem in \eqref{e:llt}.
\begin{lemma}\label{lem:tp_harnack}
For any $\varepsilon>0$, $K>0$ and $0<\delta<T<\infty$, there
exists $\varepsilon_b>0$ such that $\mathbb{P}$-a.s., for all but
finitely many $n$,
\begin{eqnarray}\label{eq:tp_harnack}
\hspace*{23pt}&&\sup_{\delta\leq t\leq T}\sup \biggl\{\frac{p_{n^2t}^\omega(nx_1,ny_1)}{p_{n^2t}^\omega(nx_2,ny_2)}{}\dvtx{} |x_i|, |y_i|
\leq K, |x_1-x_2|\leq\varepsilon_b,|y_1-y_2|\leq\varepsilon_b \biggr\}\nonumber
\\[-8pt]\\[-8pt]
&&\qquad < 1+\varepsilon.\nonumber
\end{eqnarray}
\end{lemma}
\begin{pf}
By Lemma~\ref{lem:hk_est}(g), we can assume that
the event $\{ \max_{|x|\leq Kn} U_x\leq b_n\}$ holds. So, by Lemma~\ref{lem:hk_est}(i) we get that for all $t\geq
\delta$,
\[
n^d  |p^\omega_{n^2t}(nx_1,ny_1) - p^\omega
_{n^2t}(nx_1,ny_2) | \leq C
\delta^{-(d+\theta)/2} \cdot|y_1-y_2|^\theta
\vee \bigg| \frac{b_n}{n}  \bigg|^\theta.
\]
On the other hand, by Lemma~\ref{lem:hk_est}(c), there
exists $\varepsilon_1>0$ such that for all $n$ large such that
$n^2\delta\geq b_n^2\vee n^{1+\eta} (2K)^{1+\eta}$, all
$\delta\leq t\leq T$ and $|x_1|,|y_1|\leq K$,
\[
n^d p^\omega_{n^2t}(nx_1,ny_1)\geq\varepsilon_1.
\]
Hence
\[
 \bigg|1 - \frac{p^\omega_{n^2t}(nx_1,ny_2)}{p^\omega_{n^2t}(nx_1,ny_1)} \bigg| \leq
 \frac{C\delta^{-(d+\theta)/2}}{\varepsilon_1} \cdot|y_1-y_2|^\theta\vee \bigg| \frac{b_n}{n}  \bigg|^\theta.
\]
The conclusion follows by taking $\varepsilon_b$ small enough so
that
\[
\frac{C
\delta^{-(d+\theta)/2}}{\varepsilon_1} \cdot\varepsilon_b^\theta<
\sqrt{1+\varepsilon}-1,
\]
and then interchanging the roles of $x$ and $y$ in the argument
above.
\end{pf}

\begin{pf*}{Proof of Theorem~\ref{thm:unif_llt}}
Let $\varepsilon _0>0$, to be chosen later. We first show that for any
fixed $|x|, |y|\leq K$, $\mathbb{P}$-a.s.,
\begin{eqnarray}\label{eq:tp_conv_st_end}
\frac{1}{(1+\varepsilon_0)^4}
&\leq&
\liminf_{n\rightarrow \infty}\inf_{\delta\leq t\leq T}\frac{n^dp_{n^2t}^\omega(nx,ny)}{k_{t}(x,y)}\nonumber\\[-8pt]\\[-8pt]
&\leq&\limsup_{n\rightarrow \infty}\sup_{\delta\leq t\leq T}\frac{n^dp_{n^2t}^\omega(nx,ny)}{k_{t}(x,y)}
\leq(1+\varepsilon_0)^4.\nonumber
\end{eqnarray}
The proof is similar to that in Lemma~4.2 in \cite{BC09}. First
fix an $\varepsilon_b$ so that the LHS in \eqref{eq:tp_harnack} in
Lemma~\ref{lem:tp_harnack} is bounded by ${1+\varepsilon_0}$. For
any path $\gamma\in D([0,\infty);\mathbb{R}^d)$, define the hitting
time $\sigma(\gamma)=\inf\{t{}\dvtx{} \gamma_{t}\in B(x,
\varepsilon_b)\}$. Then by the QFCLT for the VSRW $Y^{(n)}$ we get
that $\mathbb{P}$-a.s.,
\begin{eqnarray*}
&&\lim_n E_0^\omega\mathbf{1}\bigl\{Y^{(n)}_{\sigma(Y^{(n)}) + t} \in B(y,\varepsilon_b)\bigr\}\\
&&\qquad= E_0  \biggl(\mathbf{1}\{\sigma(W^*)<\infty\}\int_{z\in B(y,\varepsilon_b)} k_t\bigl(W^*_{\sigma(W^*)}, z\bigr)\,dz\biggr),
\end{eqnarray*}
where $W^*$ is the limit of the VSRW $Y^{(n)}$. So, writing
$\sigma= \sigma(Y^{(n)})$, for all large $n$,
\begin{eqnarray*}
P^0_{\omega}\bigl( Y^{(n)}_{\sigma+ t} \in B(y, \varepsilon_b) | Y^{(n)}_{\sigma}, \sigma< \infty\bigr)
 &=& \sum_{z \in B(ny, n\varepsilon _b)} p^{\omega}_{n^2t}\bigl(nY^{(n)}_\sigma, z\bigr) \\
&\ge&(1+\varepsilon _0)^{-1} |B(ny, n \varepsilon _b)|\cdot p^{\omega}_{n^2 t}\bigl(nY^{(n)}_{\sigma},ny\bigr) \\
&\ge&(1+\varepsilon _0)^{-2} |B(ny, n \varepsilon _b)|\cdot p^{\omega}_{n^2 t}(nx,ny).
\end{eqnarray*}
Note that $|B(ny, n \varepsilon _b)|\sim n^d\cdot
\operatorname{Vol}(B_\mathbb{R}(y,\varepsilon _b))$; using this and the
analogous result
for $k_t(x,y)$, we get that
\[
\limsup_n n^dp^\omega_{n^2 t}(nx, ny) \cdot P^0_{\omega}
\bigl(\sigma\bigl(Y^{(n)}\bigr)<\infty\bigr)
\leq(1+\varepsilon_0)^4 P_0\bigl(\sigma(W^*)<\infty\bigr) k_t(x,y).
\]
But by the QFCLT for the VSRW $Y^{(n)}$ again, $\lim_n P^0_{\omega}
(\sigma(Y^{(n)})<\infty) =\break P_0(\sigma(W^*)<\infty) $, hence we get
the desired upper bound. The lower bound in
\eqref{eq:tp_conv_st_end} can be proved similarly.

We now let $x,y$ vary over $B_\mathbb{R}(0,K)$. Find a finite set
$\{z_1,\ldots, z_\ell\}$ such that $B_\mathbb{R}(0,K)$ is covered by the
balls $B_\mathbb{R}(z_i,\varepsilon_b)$. By the previous argument,
$\mathbb{P}$-a.s., for all $i,j=1,\ldots,\ell$, $n^d
p_{n^2t}^\omega(nz_i,nz_j)/k_{t}(z_i,z_j)$ is bounded above by
$(1+\varepsilon _0)^4$ for all large $n$. Given $x, y\in B_\mathbb{R}(0,K)$,
choose $z_i, z_j$ so that $x\in B_\mathbb{R}(z_i,\varepsilon _b)$, $y
\in
B_\mathbb{R}(z_j,\varepsilon _b)$. Then using \eqref{eq:tp_harnack},
\begin{eqnarray*}
\frac{n^d p_{n^2t}^\omega(nx,ny)}{k_{t}(x,y)} = \frac{ n^d p_{n^2
t}^{\omega}(nz_i, nz_j) }{k_t(z_i,z_j)} \cdot\frac{ n^d p_{n^2
t}^{\omega}
(nx, ny) } { n^d p_{n^2 t}^{\omega}(nz_i, nz_j) } \cdot
\frac{k_t(z_i,z_j)}{k_t(x,y)} < (1+\varepsilon _0)^6
\end{eqnarray*}
for all large $n$. Taking $(1+\varepsilon _0)^6 < 1+\varepsilon $
gives the
upper bound in \eqref{eq:tp_unif_conv}, and the lower bound can be
proved similarly.
\end{pf*}

\subsection{Convergences after truncation}

For any given $a>0$, we introduce the following truncation of
$\mu_x$:
%
\begin{equation}\label{e:tmu_def}
\widetilde\mu _e =\widetilde\mu ^{(n)}_e= \mu_e \cdot\mathbf{1}_{\{\mu_e \leq a n^2\}},\qquad
\widetilde\mu _x=\widetilde\mu ^{(n)}_x= \sum_{y \sim x}\widetilde\mu _{xy}.
\end{equation}
Then we have
%
\begin{equation}\label{e:Etmu}
\mathbb{E}\widetilde\mu _x \sim\log(a n^2),\qquad
\mathbb{E}\widetilde\mu _x^2 \leq C an^2,
\end{equation}
where $C$ is a constant independent of $a$ and $n$. Note that
$\widetilde\mu _x$ and $\widetilde\mu _y$ are independent if $|x-y|>1$.
\begin{lemma}\label{lem:DCT_bd} Let $K>0$ and $d\geq3$.
\begin{longlist}[(\textup{a})]
\item[(\textup{a})] If $f{}\dvtx{} B_\mathbb{R}(0,K)\rightarrow \mathbb{R}$ is continuous,
then $\mathbb{P}$-a.s.,
%
\begin{equation}\label{e:sum_mu_ub}
\frac{1}{n^d \log n}\sum_{|x|\leq Kn} \widetilde\mu _x f(x/n)
\rightarrow 2\int_{ B_\mathbb{R}(0,K)} f(x)\, dx.
\end{equation}
\item[(\textup{b})] If $g{}\dvtx{} (B_\mathbb{R}(0,K))^2\rightarrow \mathbb{R}$ is
continuous, then
$\mathbb{P}$-a.s.,
%
\begin{equation}\label{e:sum_mu2_ub}
\hspace*{35pt}\frac{1}{n^{2d} (\log n)^2}\sum_{|x|,|y|\leq Kn}
\widetilde\mu _x \widetilde\mu _y  g(x/n,y/n) \rightarrow
4\int_{(B_\mathbb{R}(0,K))^2} g(x,y)\,dx\,dy.
\end{equation}
\end{longlist}
\end{lemma}

\begin{pf}
In both cases we use a straightforward
mean--variance calculation.

(a) Write $I_n$ for the LHS of \eqref{e:sum_mu_ub}. Then as
$\mathbb{E}\widetilde\mu _x \sim\log(an^2)\sim2\log n$,%
\begin{equation}\label{e:dctL1}
\hspace*{20pt}\mathbb{E}I_n= \frac{\mathbb{E}\widetilde\mu _0}{\log n} \sum_{|x|\leq Kn}
f(x/n) n^{-d} \rightarrow 2\int_{|x|\leq K} f(x)\, dx\qquad \mbox{as $n\rightarrow \infty$.}
\end{equation}
If $|x-y|\le1$, then $|\operatorname{Cov}(\widetilde\mu _x,
\widetilde\mu _y)|\leq\operatorname{Var}(\widetilde\mu _{0})$
by Cauchy--Schwarz. So
\begin{eqnarray*}
\operatorname{Var}_\mathbb{P}(I_n )
&\leq&\frac{c\|f\|_\infty^2 }{n^{2d} (\log n)^2}
\sum_{|x|\leq Kn} \operatorname{Var}(\widetilde\mu _{0}) \\
&\leq&\frac{C }{n^{d} (\log n)^2} an^2 \le\frac{C'}{n^{d-2} (\log
n)^2 }.
\end{eqnarray*}
So, for any $\varepsilon >0$ we deduce
\[
\mathbb{P}( |I_n - \mathbb{E}I_n| > \varepsilon ) \le\frac{
\operatorname{Var}_\mathbb{P}(I_n)}{\varepsilon ^2}
\le\frac{c(\varepsilon )}{n^{d-2} (\log n)^2},
\]
and so by Borel--Cantelli, we have that $|I_n - \mathbb{E}I_n| <
\varepsilon $
for all large $n$.

(b) Let $J_n$ be the left-hand side of \eqref{e:sum_mu2_ub}. Write
$B = B(0,Kn)$ and
\begin{eqnarray*}
J'_n &=& \frac{1}{n^{2d} (\log n)^2}\sum_{x,y\in B, |x-y|\le3}
\widetilde\mu _x \widetilde\mu _y g(x/n,y/n), \\
J''_n &=& \frac{1}{n^{2d} (\log n)^2}\sum_{x,y\in B, |x-y|> 3}
\widetilde\mu _x \widetilde\mu _y g(x/n,y/n).
\end{eqnarray*}
Then since $\widetilde\mu _x \widetilde\mu _y \le\widetilde\mu
_x^2 + \widetilde\mu _y^2$,
\[
\mathbb{E}|J'_n|
\le\frac{c}{n^{2d} (\log n)^2}\sum_{x\in B}\mathbb{E}\widetilde
\mu _x^2
\|g\|_\infty\le\frac{c \|g\|_\infty} {n^{d-2} (\log n)^2}.
\]
As this sum converges, by Borel--Cantelli $J'_n \rightarrow 0$
$\mathbb{P}$-a.s.

For $J''_n$ we have
\[
\mathbb{E}J''_n = \frac{ (\mathbb{E}\widetilde\mu _x)^2 }{n^{2d}
(\log n)^2}\sum_{x,y\in
B, |x-y|> 3}
g(x/n,y/n) \rightarrow 4\int_{|x|,|y|\leq K} g(x,y)\,dx\,dy.
\]
Furthermore,
%
\begin{eqnarray}\label{e:varjn}
\hspace*{30pt}\operatorname{Var}_\mathbb{P}(J''_n)&\leq&\frac{C}{n^{4d} (\log n)^4}\nonumber
\\[-8pt]\\[-8pt]
&&{}\times\sum_{ x,y \in B, |x-y|> 3}
\biggl(\sum_{x',y' \in B, |x'-y'|> 3}|\operatorname{Cov}( \widetilde\mu _x \widetilde\mu _y, \widetilde\mu _{x'}\widetilde\mu _{y'})|  \biggr).\nonumber
\end{eqnarray}
If all of $x, y, x', y'$ are at a distance greater than 1 apart in
the sum in \eqref{e:varjn}, then $ \operatorname{Cov}( \widetilde
\mu _x \widetilde\mu _y,
\widetilde\mu _{x'}\widetilde\mu _{y'})=0$. So, after relabelling,
we only have to
handle two cases: when $|x-x'|\le1 $ and $|y-y'|\le1$, and when
$|x-x'|\le1$ and $|y-y'|> 1$. Write $K'_n$ and $K''_n$ for these
two sums. Observe that in both cases, since $|x-y|>3$ and
$|x'-y'|>3$, we have $|y'-x|>1$ and $|y-x'|>1$.

In the first case,
%
\begin{equation}\label{e:cov1}
|\operatorname{Cov}( \widetilde\mu _x \widetilde\mu _y, \widetilde
\mu _{x'}\widetilde\mu _{y'})|\le
\mathbb{E}\widetilde\mu _x \widetilde\mu _{x'}\cdot\mathbb{E}\widetilde\mu _y\widetilde\mu _{y'}\le c n^4,
\end{equation}
and so
\[
K'_n \le\frac{c n^{2d} n^4}{n^{4d} (\log n)^4}\le\frac{c}{
n^{2d-4}(\log n)^4}.
\]
In the second case,
\[
|\operatorname{Cov}( \widetilde\mu _x \widetilde\mu _y, \widetilde
\mu _{x'}\widetilde\mu _{y'})|
\le\mathbb{E}\widetilde\mu _x \widetilde\mu _{x'}\cdot\mathbb{E}\widetilde\mu _y\widetilde\mu _{y'} \le c n^2 (\log n)^2,
\]
and so as the sum in $K''_n$ contains $O(n^{3d})$ terms
\[
K''_n \le\frac{c n^{3d} n^2 (\log n)^2}{n^{4d} (\log n)^4}
\le\frac{c}{ n^{d-2} (\log n)^2}.
\]
Hence $\sum_n\operatorname{Var}_\mathbb{P}(J''_n)<\infty$, proving
\eqref{e:sum_mu2_ub}.
\end{pf}

Finally we state a simple lemma which can be proved by direct
computations.
\begin{lemma}\label{lem:sum_norm} For any $K>0$,
\begin{longlist}[(\textup{a})]
\item[(\textup{a})]
\[
\sum_{1\leq|x|\leq Kn} |x|^{2-d} = O(n^2).
\]

\item[(\textup{b})]
\[
\sum_{1\leq|x|\leq Kn} |x|^{4-2d} =
\cases{O(n), &\quad\mbox{\textup{when} } $d=3$,\cr
O(\log n), &\quad\mbox{\textup{when} } $d=4$,\cr
O(1), &\quad\mbox{\textup{when} } $d\geq5$.
}
\]
\end{longlist}
\end{lemma}

\section{Estimates involving Green's functions}

For the usual simple random walk on $\mathbb{Z}^d$, $d\geq3$,
Green's function $g(x,x)$ is a positive constant for all $x$. In
our case, the best available lower bound [see Lemma~\ref{lem:hk_est}(e)] gives that $\mathbb{P}$-a.s.,
for all
large $n$, and for all $|x|\leq Kn$, $g^{\omega}(x,x)\geq C/(\log
n)^{(d-2)/\eta}$. As this is not quite strong enough for the
truncation arguments in the next section, we now derive some more
precise bounds on sums of Green's functions in a ball.

Recall that $E_d$ denotes the set of edges in $\mathbb{Z}^d$, and in
Lemma~\ref{lem:hk_est}(g) we defined $b_n = c_{13} (\log
n)^{1/\eta}$. For $e=\{x_e, y_e\}\in E_d$, let
$B(e,r)=B(x_e,r)\cap B(y_e,r)$.
For $e=\{x_e,y_e\} \in E_d$ and $z \in\mathbb{Z}^d$, let
%
\begin{eqnarray}
\gamma_n(e)&=&C_{\mathrm{eff}}[\{x_e,y_e\},B(e,b_n)^c],\\
\gamma_n(z)&=&C_{\mathrm{eff}}[z,B(z,b_n+1)^c],
\end{eqnarray}
where $C_{\mathrm{eff}}[A,B]$ denotes the effective conductivity
between the
sets $A$ and $B$ (see (3.8) in \cite{BC09} or \cite{LPW}, Section~9.4). Note that both $\gamma _n(e)$ and $\gamma _n(x)$ are
decreasing in $n$, and $\gamma _\infty(e):=\lim_n \gamma _n(e)$ is the
effective conductivity between $e$ and infinity while
$\gamma _\infty(x):=\lim_n \gamma _n(x)$ is equal to $1/g^{\omega
}(x,x)$. By
\cite{BC09}, Lemma~6.2, for any $k\geq1$, $\lim_n \mathbb{E}
\gamma _n(e)^k<~\infty$. Note further that $\mu_e$ and $\gamma_n(e)$
are independent, and also that $\gamma_n(e)$ and $\gamma_n(e')$
are independent if $|e-e'| \ge2b_n+1$. When $d\geq3$, by Lemma~\ref{lem:hk_est}(e), $g^{\omega}(x,x)<C<\infty$, and
hence
%
\begin{equation}\label{eq:gm_n_lbd}
\gamma_n(e)\geq\gamma_n(x)\geq
\gamma_\infty(x)=1/g^{\omega}(x,x)\geq1/C>0.
\end{equation}
Let $a_p$ be large enough so that $\mathbb{P}(\mu_e > a_p) < p_c(d)$
where $p_c(d)$ is the critical probability for bond percolation in
$\mathbb{Z}^d$. Let $\mathcal{C}(e)$ denote the cluster
containing~$e$ in the
bond percolation process for which $\{ e \mbox{ is open} \} =\{
\mu_e> a_p\}$. Then we have (see \cite{Grim}, Theorems 6.75 and 5.4)
%
\begin{eqnarray}\label{eq:sub_perc_bd}
\mathbb{P}( |\mathcal{C}(e)| > m ) &\le&\exp( - c_1 m ),\nonumber
\\[-8pt]\\[-8pt]
\mathbb{P}\bigl( \operatorname{diam}(\mathcal{C}(e)) > m\bigr) &\le&\exp( -c_2 m),\qquad\mbox{for all } m\ge 1\nonumber
\end{eqnarray}
Let
\begin{eqnarray*}
F_n(e) = \bigl\{ \operatorname{diam}(\mathcal{C}(e)) \ge\fract12b_n\bigr\}, \qquad
\gamma _n'(e) = \gamma _n(e)\cdot\mathbf{1}_{F_n(e)^c}.
\end{eqnarray*}
\begin{lemma}\label{lem:Ceff-exp}
\textup{(a)} For any $K>0$, $\mathbb{P}$-a.s., for all sufficiently large $n$,
$\gamma _n(e)=\gamma '_n(e)$ for all $e \in B(0,2Kn)$.
\begin{longlist}[(\textup{b})]
\item[(\textup{b})] There exists $\theta>0$ and $\Gamma=\Gamma(\theta)<\infty$
such that for all $n$,
\[
\mathbb{E}e^{\theta\gamma '_n(e)} < \Gamma.
\]
\item[(\textup{c})] There exists $C=C(d)>0$ such that for any $K>0$,
$\mathbb{P}$-a.s., for all large~$n$,
\[
\inf_{|x|\leq Kn} g^{\omega}(x,x)\geq C/\log n.
\]
\end{longlist}
\end{lemma}

\begin{pf}
(a) First note that
%
\begin{equation}\label{e:Fnbnd}
\mathbb{P}\biggl( \bigcup_{e \in B(0,2Kn) } F_n(e)\biggr) \le c n^d \exp(-c_2b_n/2)
= c \exp\bigl( d \log n - c' (\log n)^{1/\eta}\bigr).\hspace*{-30pt}
\end{equation}
Since $\eta<1$ the RHS in \eqref{e:Fnbnd} is summable, so that,
for all but finitely many~$n$, $\gamma _n(e)=\gamma '_n(e)$ for all $ e \in B(0, 2 Kn)$.

(b) On $F_n(e)^c$ the cluster $\mathcal{C}(e)$ is contained in $B(e,
b_n)$, and each bond from $\mathcal{C}(e)$ to $\mathcal{C}(e)^c$ has
conductivity
less than $a_p$. Since there are at most $2 d |\mathcal{C}(e)|$ such
bonds, we deduce that $\gamma_n(e) \le d a_p |\mathcal{C}(e)|$. So,
%
\begin{equation}\label{e:gam-pt}
\mathbb{P}\bigl( \gamma '_n(e)> {\lambda}\bigr)
\le\mathbb{P}\bigl( d a_p|\mathcal{C}(e)|>{\lambda}\bigr)\le\exp( -c {\lambda}).
\end{equation}

(c)
Using \eqref{e:gam-pt} it is enough to consider
\begin{eqnarray*}
\mathbb{P} \Bigl( \max_{e \in B(0, Kn)} \gamma '_n(e) > {\lambda}\log n \Bigr)\le c' n^d e^{- c {\lambda}\log n}
\end{eqnarray*}
which is summable when ${\lambda}$ is large enough.
\end{pf}

For any $0<a<b\leq\infty$, define the sets
%
\begin{equation}
E_n(a,b)=\{ e{}\dvtx{} an^2 \le\mu_e < b n^2 \}.
\end{equation}

Let $m_n$ be chosen later with $m_n \ge3 b_n$. We tile $\mathbb{Z}^d$
with cubes of the form $Q= [0,m_n-1]^d + m_n \mathbb{Z}^d$ so that each
cube contains $m_n^d$ vertices. Let $z_i$, $1\le i \le d$, be the
unit vectors in $\mathbb{Z}^d$, and given a cube $Q$ in the tiling let
\[
E(Q)= \bigl\{ \{x,x+z_i\}, x \in Q, 1\le i \le d  \bigr\};
\]
it is clear that $E(Q)$ gives a tiling of $E_d$, and that
$|E(Q)|=d m_n^d$ for each $Q$. Let $K>0$ be fixed, and let $\mathcal{Q}_n$
be the set of $Q$ such that $Q \cap B(0,Kn+1) \neq\varnothing$. We
have $|\mathcal{Q}_n| \asymp(Kn/m_n)^d$.
\begin{lemma}[(See \cite{BC09}, Lemma~6.3)]\label{l:homog}
Let $a, K, \delta>0$ be fixed.
\begin{longlist}[\textup{(a)}]
\item[\textup{(a)}] Suppose that $Kn/\sqrt{d}\geq m_n \ge n^{\theta_1}$ for some
$\theta_1>2/d$. Then there exists ${\lambda}>0$ such that $\mathbb{P}$-a.s.,
for all but finitely many $n$,
%
\begin{equation}\label{e:Emeanbound}
\max_{ Q \in\mathcal{Q}_n} \sum_{e \in E(Q)\cap E_n(a,\infty) }
\gamma_n(e)
\le{\lambda}m_n^d (an^2)^{-1}\mathbb{E}\gamma _n(e).
\end{equation}
\item[\textup{(b)}] Let $\theta_2 < 1/d$. Then $\mathbb{P}$-a.s., $B(0, n^{\theta_2})
\cap E_n(a,\infty) =\varnothing$ for all but finitely many $n$.
\end{longlist}
\end{lemma}

\begin{pf} (a)
By Lemma~\ref{lem:Ceff-exp}(a) it is enough to bound the sum
\eqref{e:Emeanbound} with $\gamma '_n(e)$ instead of $\gamma _n(e)$. Let
$Q \in\mathcal{Q}_n$. We divide $E(Q)$ into disjoint sets $(E(Q,j), j
\in
J)$ such that if $e$ and $e'$ are distinct edges in $E(Q,j)$, then
$|e-e'| \ge3 b_n-2$, each $|E(Q,j)| = (m_n/3b_n)^d:=N_n$, and
$|J| \sim d (3b_n)^d$.

Let $\eta_e = \mathbf{1}_{ (\mu_e > an^2) }$, $p_n= \mathbb{E}\eta_e
\sim1/(2d)\cdot1/(a n^2)$, and
\[
\xi_j = \sum_{e \in E(Q,j) } \gamma '_n(e) \eta_e.
\]
Then the r.v. $( \gamma '_n(e), \eta_e, e \in E(Q,j))$ are
independent, and so if $\theta$ and $\Gamma$ are as in
Lemma~\ref{lem:Ceff-exp},
\[
\mathbb{E}e^{\theta\xi_j} \le\bigl(1 + p_n (\Gamma-1)\bigr)^{N_n}\leq e^{N_n p_n
(\Gamma-1) }.
\]
Hence for any ${\lambda}>0$, writing $ \mathbb{E}\xi_j = N_n p_n
\mathbb{E}
\gamma '_n(e)$,
\begin{eqnarray*}
\mathbb{P}( \xi_j > {\lambda}\mathbb{E}\xi_j )
&\leq&\exp\bigl(-{\lambda}\theta N_n p_n \mathbb{E}\gamma'_n(e) + N_n p_n (\Gamma-1)\bigr)\\
&=&\exp \bigl(-N_n p_n  \bigl({\lambda}\theta\mathbb{E}\gamma'_n(e)- \Gamma+1 \bigr) \bigr).
\end{eqnarray*}
By \eqref{eq:gm_n_lbd},
\[
\mathbb{E}\gamma '_n(e)\geq1/C\cdot\mathbb{P}(
F_n(e)^c)\rightarrow 1/C,
\]
hence there exists ${\lambda}>0$ such that for all $n$ large,
${\lambda}\theta\mathbb{E}\gamma '_n(e)- \Gamma+1\geq1$, and so
\begin{eqnarray*}
\mathbb{P} ( \xi_j > {\lambda}\mathbb{E}\xi_j  ) \le
e^{-N_n p_n}.
\end{eqnarray*}
Thus
\[
\mathbb{P} \biggl( \sum_{j \in J} \xi_j > {\lambda}m_n^d p_n \mathbb{E}\gamma '_n(e) \biggr)
\le d (3 b_n)^d e^{-N_n p_n},
\]
and so since $|\mathcal{Q}_n| \le cn^d$ and $N_np_n\geq n^\varepsilon
$ for some
$\varepsilon >0$, \eqref{e:Emeanbound} follows by Borel--Cantelli.

(b) We have
\[
\mathbb{P}\bigl( B(0, n^{\theta_2}) \cap E_n(a,\infty) \neq\varnothing\bigr)
\le c n^{d \theta_2} (an^2)^{-1} \leq c n^{d \theta_2 -2};
\]
so again the result follows using Borel--Cantelli.
\end{pf}

\section[Proof of Theorem 2]{Proof of Theorem~\protect\ref{thm:conv}}
\begin{lemma}\label{lem:skorohod_marginal}
Let ${\omega}\in\Omega$. If for each $t\geq0$,
%
\begin{equation}\label{eq:conv_clock_marginal}
S^{(n)}_t \rightarrow 2t\qquad  \mbox{in $P^0_{\omega}$-probability},
\end{equation}
then \eqref{eq:conv_clock} holds.
\end{lemma}
\begin{pf}
Note that the LHS and RHS are both increasing processes, and the
RHS is continuous and deterministic. The conclusion then follows
from Theorem~VI.3.37 in \cite{JS03}.
\end{pf}
\begin{lemma}\label{lem:conv_trunc}
For each $\varepsilon>0$ and $T>0$, there exist $K>0$ and $a>0$
such that for $\mathbb{P}$-a.a. ${\omega}$, for all $t\leq
T$, the
following two inequalities hold:
%
\begin{eqnarray}\label{ineq:big_ball}
\limsup_n P^0_{\omega} \biggl(\frac{1}{n^2\log n}\sum_{|x|\geq Kn }\int_0^{n^2 t} \mu_x\cdot\mathbf{1}_{\{Y_s =x\}}\,ds >0  \biggr)
&\leq&\varepsilon;\hspace*{-40pt}
\\\label{ineq:big_weight}
\qquad \limsup_n P^0_{\omega} \biggl( \frac{1}{n^2\log n}\sum_{|x|\leq Kn }
\int_0^{n^2 t} \mu_x\cdot
\mathbf{1}_{\{\mu_x\geq a n^2\}}
\mathbf{1}_{\{Y_s=x\}}\,ds >0  \biggr)
&\leq&\varepsilon.
\end{eqnarray}
\end{lemma}

\begin{pf}
Write $F_K$ for the event in \eqref{ineq:big_ball}. Then by Lemma~\ref{lem:hk_est}(d),
\[
P^0_{\omega}(F_K) \le P^0_{\omega}\bigl( \tau(0,Kn) < n^2t\bigr) \le c_8 \exp
(-c_9 K^2/t),
\]
provided that $Kn > U_0$. So, taking $K$ sufficiently large,
\eqref{ineq:big_ball} holds for all sufficiently large $n$.

Choose $\theta_1 = (2+\varepsilon _1)/d$, $\theta_2 =
(1-\varepsilon _2)/(d-2)$
where $\varepsilon _1>0, \varepsilon _2>2/d$ (so that $\theta
_2<1/d$) and $\varepsilon _1
+ \varepsilon _2<1$. Let $m_n = n^{\theta_1}$, and $\mathcal{Q}_n$
be as in Lemma~\ref{l:homog}. Let $n$ be large enough so that
\eqref{e:Emeanbound} holds, and also that
%
\begin{equation}\label{eq:smball_reg}
B(0,n^{\theta_2}) \cap E_n(a,\infty)=\varnothing.
\end{equation}
Then
%
\begin{equation}
P^0_{\omega}\bigl(Y \mbox{ hits } E_n (a,\infty) \cap B(0,Kn) \bigr)
\le\sum_{Q \in\mathcal{Q}_n} \sum_{x \in E_n(a,\infty) \cap Q}
\frac{
g^{\omega}(0,x)}{g^{\omega}(x,x)}.
\end{equation}
For $x \in E_n(a,\infty)$, if $e_x$ is an edge containing $x$,
then by \eqref{eq:gm_n_lbd} $1/g^{\omega}(x,x) \le\gamma_n(e_x)$. By
\eqref{eq:smball_reg} and \eqref{eq:Gn_gen_bd} we can bound
$g^{\omega}(0,x)$ by $c |x|^{2-d}$.

Let $\mathcal{Q}_n'$ be the set of $Q \in\mathcal{Q}_n$ such that
$|x| \ge m_n/2$
for all $x \in Q$. Let first $Q \in\mathcal{Q}_n \backslash\mathcal{Q}_n'$. Then
by Lemma~\ref{l:homog} and \eqref{eq:smball_reg},
\begin{eqnarray*}
\sum_{x \in E_n (a,\infty) \cap Q} \frac{ g^{\omega
}(0,x)}{g^{\omega}(x,x)}
&\le&\max_{x \in E_n(a,\infty) \cap Q} c |x|^{2-d}
\sum_{x \in E_n (a,\infty) \cap Q} \gamma _n(e_x) \\
&\le& C n^{\theta_2 (2-d)}\cdot{\lambda}m_n^d (a n^2)^{-1}
\le C' n^{\varepsilon _1 + \varepsilon _2 -1}.
\end{eqnarray*}
So, since there are only $2^d$ cubes in $\mathcal{Q}_n -\mathcal{Q}_n'$
and $\varepsilon_1+\varepsilon_2<1$ by the choices of $\varepsilon_1$
and $\varepsilon_2$,
%
\begin{equation}\label{e:q-qp}
\lim_n \sum_{Q \in\mathcal{Q}_n-\mathcal{Q}_n'} \sum_{x \in E_n
(a,\infty) \cap Q}
\frac{ g^{\omega}(0,x)}{g^{\omega}(x,x)} =0.
\end{equation}

Now let $Q \in\mathcal{Q}_n'$, and let $x_Q$ be the point in $Q$ closest
to 0. Then if $Q \in\mathcal{Q}_n'$,
%
\begin{eqnarray}\nonumber
\sum_{x \in E_n (a,\infty) \cap Q} \frac{ g^{\omega
}(0,x)}{g^{\omega}(x,x)}
&\le& c \sum_{x \in E_n (a, \infty) \cap Q} |x|^{2-d} \gamma _n(e_x)
\\
&\le& c |x_Q|^{2-d}\cdot{\lambda}m_n^d (a n^2)^{-1}
\\
&\le& c' {\lambda}a^{-1} n^{-2} \sum_{x \in Q} |x|^{2-d}.\nonumber
\end{eqnarray}
So, summing over $Q \in\mathcal{Q}_n'$,
\begin{eqnarray*}
P^0_{\omega}\biggl(Y \mbox{ hits } E_n (a,\infty) \cap\biggl(\bigcup_{Q \in\mathcal{Q}_n'}Q\biggr) \biggr)
&\le& c{\lambda}a^{-1} n^{-2} \sum_{x \in B(0,(K+1)n )} (1 \vee|x|)^{2-d}
\\
&\le& c' {\lambda}(K+1)^2   a^{-1},
\end{eqnarray*}
and so taking $a$ large enough and noting \eqref{e:q-qp},
\eqref{ineq:big_weight} follows.
\end{pf}

By Lemma~\ref{lem:conv_trunc} to prove \eqref{eq:conv_clock} it
suffices to consider the convergence of
%
\begin{eqnarray}\label{eq:main_term}
\widetilde S^{(n)}_t &=& \frac{1}{n^2\log n} \sum_{|x|\leq Kn}
\widetilde\mu _x \cdot\int_0^{n^2 t}\mathbf{1}_{\{Y_s=x\}}\,ds\nonumber
\\[-8pt]\\[-8pt]
&=&\frac{1}{\log n}\sum_{|x|\leq Kn}\widetilde\mu _x \cdot\int_0^{t}
\mathbf{1}_{\{Y_{n^2s}=x\}}\,ds,\nonumber
\end{eqnarray}
where $\widetilde\mu _x$ is as in \eqref{e:tmu_def}. Taking expectations
with respect to $P^0_{\omega}$ we have
%
\begin{eqnarray}\label{e:sndef}
\nonumber
E^0_{\omega}\widetilde S^{(n)}_t
 &=& \frac{1}{n^2 \log n}\sum_{|x|\le Kn} \widetilde\mu _x
\cdot\int_0^{n^2 t} p^\omega_{s}(0,x)\,ds \nonumber\\[-8pt]\\[-8pt]
 &=& \frac{1}{\log n}\sum_{|x|\leq
Kn} \widetilde\mu _x
\cdot\int_0^{t} p^\omega_{n^2 r}(0,x)\,dr.\nonumber
\end{eqnarray}
\begin{lemma}\label{lem:ini_ignore}
For any $\varepsilon>0$, there exists $\delta>0$ such that,
$\mathbb{P}$-a.s. for all sufficiently large~$n$,
%
\begin{equation}\label{e:initsn}
E^0_{\omega}\widetilde S^{(n)}_\delta\leq\varepsilon.
\end{equation}
\end{lemma}

\begin{pf} By Lemma~\ref{lem:hk_est}(g),
we can assume $n$ is large enough so that \break
$\{\max_{|x|\leq Kn} U_x\leq b_n\}$.
Hence, by Lemma~\ref{lem:hk_est}(b), if
$|x|\vee\sqrt{t}\geq b_n$, then
\begin{eqnarray*}
p_t^\omega(0,x) \leq
 \cases{ c_4 t^{-d/2}\exp(-c_5 |x|^2/t), &\quad \mbox{when }$t\geq|x|$,\cr
c_4\exp(-c_5 |x|), &\quad\mbox{when }$t\leq|x|$.}
\end{eqnarray*}
Hence, by decomposing according to whether $|x| < b_n$ or $|x|\ge
b_n$, we obtain
%
\begin{eqnarray}\nonumber
E^0_{\omega}\widetilde S^{(n)}_\delta
&=&\frac{1}{n^2 \log n}\sum_{|x|\le Kn} \widetilde\mu _x\cdot\int_0^{n^2 \delta}
p^\omega_{ s}(0,x)\,ds \\\label{e:snsum-3}
&\le&\frac{1}{n^2 \log n}\sum_{|x|\le b_n} \widetilde\mu _x \cdot
\int_0^{n^2\delta} c (1\vee s )^{-d/2}\,ds\\\label{e:snsum-1}
&&{}+ \frac{1}{n^2 \log n} \sum_{b_n\leq|x|\leq Kn}
\widetilde\mu _x \int_0^{|x|} c_4 e^{-c_5|x|}\,ds\\\label{e:snsum-2}
&&{}+\frac{1}{n^2 \log n}\sum_{b_n\leq|x|\le Kn}
\widetilde\mu _x\cdot\int_{|x|}^{n^2 \delta} c_4
s^{-d/2} e^{-c_5|x|^2/s }\,ds.
\end{eqnarray}
Write $\xi^{(i)}_n$, $i=1,2, 3$, for the terms in
\eqref{e:snsum-3}--\eqref{e:snsum-2}. Since the integral in
\eqref{e:snsum-3} is bounded by $\int_0^{\infty} c (1\vee s)^{-d/2}\, ds<\infty$, we have
\[
\mathbb{E}\xi^{(1)}_n \le c \frac{b_n^{d}}{n^2 \log n} \mathbb{E}\widetilde\mu _x
\le c n^{-2} (\log n)^{d/\eta}.
\]
Similarly for \eqref{e:snsum-1} we have
\[
\mathbb{E}\xi^{(2)} \le c n^{-2} \sum_{|x| \le Kn} c_4|x| e^{-c_5|x|}
\le c' n^{-2}.
\]
As these sums converge it follows from Borel--Cantelli that $
\xi^{(i)}_n \le\varepsilon /3$ for all large~$n$, for $i=1,2$.

It remains to control \eqref{e:snsum-2}. First note that when $s
\ge1$,
%
\begin{equation}\label{eq:tp_normal}
\sum_{ x \in\mathbb{Z}^d} s^{-d/2}e^{- \kappa|x|^2/s}
\le C(\kappa).
\end{equation}
So, interchanging the order of the sum and integral in
\eqref{e:snsum-2},
\[
\mathbb{E}\xi^{(3)}_n\leq\frac{C}{n^2\log n} \mathbb{E}\widetilde
\mu _0 \cdot n^2\delta
\le C' \delta.
\]
Setting $t=s/|x|^2$ we have
%
\begin{eqnarray}\label{eq:normal_bd}
\hspace*{25pt}\int_{|x|}^{n^2 \delta} c_4 s^{-d/2}e^{-c_5|x|^2/s }\,ds
\le C |x|^{2-d} \int_0^\infty t^{-d/2 } e^{-c_5/t}\,dt
\le C |x|^{2-d}.
\end{eqnarray}
Hence, applying Lemma~\ref{lem:sum_norm} we get
\[
\operatorname{Var}_\mathbb{P} \bigl(\xi^{(3)}_n \bigr)
\leq\frac{C }{n^ 4 (\log n)^2}\cdot
\sum_{b_n \le|x|\leq Kn} an^2 |x|^{4-2d}
\le\frac{C}{n(\log n)^2 }.
\]
By Chebyshev's inequality and Borel--Cantelli we then get that for
$\delta$ small enough, $\mathbb{P}$-a.s. for all sufficiently
large~$n$, $\xi^{(3)}_n\leq\varepsilon/3$.
\end{pf}
\begin{prop}\label{prop:conv_exp}
Let
%
\begin{equation}\label{e:A1_def}
A_1(K,t,\delta)= \int_{|y|\leq K} \int_\delta^t k_s(x)\,dx\,ds.
\end{equation}
When $d\geq3$, for any $K>0$, $0<\delta<T<\infty$, and
$t\in(\delta, T]$, $\mathbb{P}$-a.s.,
%
\begin{equation}\label{eq:bd_exp}
\lim_{n\rightarrow \infty} E^0_{\omega} \bigl(\widetilde
S^{(n)}_t-\widetilde S^{(n)}_\delta \bigr)=2
A_1(K,t,\delta).
\end{equation}
\end{prop}
\begin{pf}
By Lemma~\ref{lem:DCT_bd}(a), it suffices to show that
$\mathbb{P}$-a.s.,
\[
\frac{1}{\log n}\sum_{|x|\leq
Kn} \widetilde\mu _x
\cdot\int_\delta^{t}  \bigl(p^\omega_{n^2 s}(0,x) - n^{-d}
k_s(x/n) \bigr)\,ds \rightarrow 0.
\]
The LHS is bounded in absolute value by
\[
\frac{1}{n^d\log n}\sum_{|x|\leq
Kn} \widetilde\mu _x\cdot T\sup_{x\in\mathbb{Z}^d}\sup_{s\geq
\delta}
|n^d p^\omega_{n^2 s}(0,x) - k_s(x/n)|.
\]
This converges to 0 $\mathbb{P}$-a.s. by Lemmas~\ref{lem:DCT_bd}(a) and \ref{lem:hk_est}(h).
\end{pf}
\begin{prop}\label{prop:conv_var_qu}
When $d\geq3$, for any $\varepsilon>0$, $K>0$,
$0<\delta<T<\infty$, and
$t\in(\delta, T]$, $\mathbb{P}$-a.s.,
%
\begin{eqnarray}\label{e:sn2_conv}
&&\limsup_n E^0_{\omega} \bigl(\widetilde S^{(n)}_t-\widetilde S^{(n)}_\delta \bigr)^2\nonumber\\[-8pt]\\[-8pt]
&&\qquad\leq \varepsilon+ 8(1+\varepsilon) \int_{|x|,|y|\leq K}
\int_\delta^t k_s(x) \int_0^{t-s} k_r(x,y)\, dr\,  ds\,  dx\,  dy.\nonumber
\end{eqnarray}
\end{prop}

\begin{pf}
Using the Markov property and the symmetry of $Y$,
\begin{eqnarray*}
&&E^0_{\omega}\bigl(S^{(n)}_t - S^{(n)}_\delta\bigr)^2\\
&&\qquad= \frac{2}{(\log n)^2}
\biggl(\sum_{|x|,|y|\leq Kn}\widetilde\mu _x\widetilde\mu _y\cdot\int_\delta^{t} p^\omega
_{n^2 s}(0,x)
\int_0^{t-s}p^\omega_{n^2 r}(x,y)\, dr\,  ds  \biggr).
\end{eqnarray*}

We begin by proving that, given $\varepsilon>0$, there exists
$\delta_1>0$ such that $\mathbb{P}$-a.s., for all large~$n$,
%
\begin{equation}\label{e:claim}
\frac{2}{(\log n)^2}\sum_{|x|,|y|\leq
Kn}\widetilde\mu _x\widetilde\mu _y\cdot\int_\delta^{t} p^\omega_{n^2 s}(0,x)
\int_0^{\delta_1}p^\omega_{n^2 r}(x,y)\,dr\, ds \leq\varepsilon.
\end{equation}

By Lemma~\ref{lem:hk_est}(a) we have $p^{\omega}_{n^2s}(0,x) \le c
n^{-d} $ for all $s\geq\delta$ and so the LHS of \eqref{e:claim}
is bounded by
%
\begin{eqnarray}
&&\frac{C}{n^d (\log n)^2}  \sum_{|x|,|y|\leq Kn}\widetilde\mu _x\widetilde\mu _y\int_0^{\delta_1}p^\omega
_{n^2 r}(x,y)\,dr \\\label{e:vsum1}
&&\qquad= \frac{C}{n^{d+2} (\log n)^2}\sum_{|x|,|y|\leq
Kn, |x-y|>1}  \widetilde\mu _x \widetilde\mu _y \int_0^{n^2
\delta_1} p^\omega_r
(x,y)\,dr\\\label{e:vsum2}
&&\qquad\quad{}+ \frac{C}{n^{d+2} (\log n)^2}\sum_{|x|,|y|\le Kn, |x-y|\le1}
\widetilde\mu _x\widetilde\mu _y \int_0^{n^2 \delta_1} p^\omega
_r(x,y)\,  dr.
\end{eqnarray}
Write $A_n$ and $B_n$ for the terms in \eqref{e:vsum1} and
\eqref{e:vsum2}.

The first term can be handled in the same way as in Lemma~\ref{lem:ini_ignore}. Let $B=B(0,Kn)$, and write $A_n= A_n^{(1)} +
A_n^{(2)} + A_n^{(3)}$ where
%
\begin{eqnarray}\label{e:asum1}
A_n^{(1)} &=& \frac{C}{n^{d+2} (\log n)^2}
\sum_{x,y \in B,  1< |x-y|\le b_n} \widetilde\mu _x \widetilde\mu _y
\int_0^{n^2 \delta_1} p^\omega_r(x,y)\,dr, \\\label{e:asum2}
A_n^{(2)} &=& \frac{C}{n^{d+2} (\log n)^2}
\sum_{ x,y \in B,  |x-y|\ge b_n} \widetilde\mu _x \widetilde\mu _y\int_0^{|x-y| } p^\omega_r(x,y)\,dr, \\\label{e:asum3}
A_n^{(3)} &=& \frac{C}{n^{d+2} (\log n)^2}
\sum_{x,y \in B,   |x-y|\ge b_n} \widetilde\mu _x \widetilde\mu _y
\int_{|x-y|}^{n^2 \delta_1} p^\omega_r(x,y)\,dr.
\end{eqnarray}

For \eqref{e:asum1} we have
\begin{eqnarray*}
\mathbb{E}A^{(1)}_n &\le&\frac{C}{n^{d+2} (\log n)^2}
\sum_{x,y \in B, 1< |x-y|< b_n} \mathbb{E}(\widetilde\mu _x
\widetilde\mu _y)
\int_0^\infty c_4 (1\vee s)^{-d/2}\,ds \\
&\le&\frac{C}{n^{d+2} } K^d n^d b_n^d
\\
&\le& c \frac{ (\log n)^{d/\eta}}{n^2},
\end{eqnarray*}
and since this sum converges, we have $A^{(1)}_n \le\varepsilon /4$ for
all large $n$, $\mathbb{P}$-a.s. The term $ \mathbb{E}A^{(2)}_n$ is
bounded in
the same way as was the term $\xi^{(2)}_n$ in Lemma~\ref{lem:ini_ignore}.

For \eqref{e:asum3},
%
\begin{eqnarray}\label{e:a31}
A^{(3)}_n &\le&\frac{C}{n^{d+2} (\log n)^2}\nonumber
\\[-8pt]\\[-8pt]
&&{}\times\sum_{x,y \in B, |x-y|>b_n} \widetilde\mu _x\widetilde\mu _y \int_{|x-y|}^{n^2 \delta_1}
c_4 s^{-d/2} \exp(-c_5|x-y|^2/s)\, ds.\hspace*{-30pt}\nonumber
\end{eqnarray}
Using \eqref{eq:tp_normal} we have
\[
\mathbb{E}A^{(3)}_n \le\frac{C }{n^{d+2} (\log n)^2} \cdot n^d
(\mathbb{E}\widetilde\mu _{0})^2\cdot n^2\delta_1 = O(\delta_1).
\]
We now bound $\operatorname{Var}_\mathbb{P}(A^{(3)}_n)$. By \eqref
{eq:normal_bd}, the
integral in \eqref{e:a31} is bounded by $c|x-y|^{2-d}$, so
\begin{eqnarray*}
\operatorname{Var}_\mathbb{P}\bigl(A^{(3)}_n\bigr)
&\le&\frac{C} {n^{2d+4 } (\log n)^4}  \\
&&{}\times\sum_{x_1, y_1 \in B,|x_1-y_1|>b_n}\sum_{x_2, y_2 \in B, |x_2-y_2|>b_n}
|x_1-y_1|^{2-d} |x_2-y_2|^{2-d}
\\
&&\hspace*{164pt}{}\times|\operatorname{Cov}(\widetilde\mu _{x_1} \widetilde\mu _{y_1},
\widetilde\mu _{x_2}\widetilde\mu _{y_2})|.
\end{eqnarray*}
We now bound this sum in the same way as was done for the variance
in Lemma~\ref{lem:DCT_bd}(b). Let
\begin{eqnarray*}
\mathcal{C}_1 &=& \{ (x_1,x_2,y_1,y_2) \in B^4{}\dvtx{}
|x_i-y_i|>b_n, i=1,2, |x_1-x_2|\le1, |y_1-y_2|\le1 \},\\
\mathcal{C}_2 &=&\{ (x_1,x_2,y_1,y_2) \in B^4{}\dvtx{}
|x_i-y_i|>b_n, i=1,2, |x_1-x_2|\le1, |y_1-y_2| > 1 \}.
\end{eqnarray*}
Note that if $|x_1-x_2|\le1$, then since $|x_i-y_i|>b_n$, none
of the $y_i$ can be within distance~1 of $x_j$.
If $(x_1, \ldots, y_2) \in\mathcal{C}_1$, then $|\operatorname
{Cov}(\widetilde\mu _{x_1} \widetilde\mu _{y_1},
\widetilde\mu _{x_2}\widetilde\mu _{y_2})| \le c n^4$, \vspace*{2pt}while if
$(x_1, \ldots, y_2) \in
\mathcal{C}_2$, then $|\operatorname{Cov}(\widetilde\mu _{x_1}
\widetilde\mu _{y_1}, \widetilde\mu _{x_2}\widetilde\mu _{y_2})|
\le c(\log n)^2 n^2$. So,
\begin{eqnarray*}
&&\frac{C}{n^{2d+4 } (\log n)^4}  \sum_{ (x_1, \ldots, y_2) \in \mathcal{C}_1}
|x_1-y_1|^{2-d} |x_2-y_2|^{2-d}\cdot
|\operatorname{Cov}(\widetilde\mu _{x_1} \widetilde\mu _{y_1},
\widetilde\mu _{x_2}\widetilde\mu _{y_2})| \\
&&\qquad\le\frac{C} {n^{2d+4 } (\log n)^4}
\sum_{ x_1, y_1 \in B} (1\vee|x_1-y_1|)^{4-2d} c n^4 \\
&&\qquad\le\frac{C} {n^{2d } (\log n)^4} n^d
\max_{x_1 \in B} \sum_{y_1 \in B(x, 2Kn)} (1\vee|x_1-y_1|)^{4-2d} \\
&&\qquad\le\frac{C n } {n^{d } (\log n)^4},
\end{eqnarray*}
where in the last inequality we used Lemma~\ref{lem:sum_norm}(b).

Also,
\begin{eqnarray*}
&&\frac{C} {n^{2d+4 }(\log n)^4} \sum_{ (x_1, \ldots, y_2) \in\mathcal{C}_2}
|x_1-y_1|^{2-d} |x_2-y_2|^{2-d} |\operatorname{Cov}(\widetilde\mu
_{x_1} \widetilde\mu _{y_1}, \widetilde\mu
_{x_2}\widetilde\mu _{y_2})| \\
&&\qquad \le\frac{C} {n^{2d+2} (\log n)^2} \sum_{ (x_1, \ldots, y_2) \in
\mathcal{C}_2}
|x_1-y_1|^{2-d} |x_2-y_2|^{2-d} \\
&&\qquad\le\frac{C}{n^{2d+2} (\log n)^2} \sum_{ x_1 \in B}\sum_{y_1, y_2 \in B(x, 2Kn)}
(1\vee|x_1-y_1|)^{2-d} (1\vee|x_1-y_2|)^{2-d} \\
&&\qquad\le\frac{C} {n^{d+2} (\log n)^2} \biggl( \sum_{y_1 \in B(0, 2Kn)} (1\vee|y_1|)^{2-d}  \biggr)^2 \\
&&\qquad\le\frac{Cn^4 } {n^{d+2} (\log n)^2} = \frac{C}{n^{d-2} (\log n)^2}.
\end{eqnarray*}
Thus $\sum_n \operatorname{Var}_\mathbb{P}(A^{(3)}_n) < \infty$,
and so if $\delta_1$
is small enough then by Chebyshev's inequality and Borel--Cantelli,
$\mathbb{P}$-a.s. for all sufficiently large $n$,
$A^{(3)}_n\leq~\varepsilon/4$.

To finish the proof of \eqref{e:claim}, it remains to bound the
term \eqref{e:vsum2}. By Lemma~\ref{lem:hk_est}(a), $ \int_0^{n^2
\delta_1}p^\omega_{r}(x,y)\,  dr \leq C$. Therefore by
Cauchy--Schwarz,
\begin{eqnarray*}
B_n&=& \frac{C}{n^{d+2} (\log n)^2}\sum_{|x|\leq
Kn, |y-x|\leq1} \widetilde\mu _x \widetilde\mu _y
\int_0^{n^2\delta_1}p^\omega_{n^2r}(x,y)\,dr\\
&\leq&\frac{C}{n^{d+2} (\log n)^2}\sum_{|x|\leq Kn} \widetilde\mu _x^2.
\end{eqnarray*}
Hence
\[
\mathbb{E}B_n\leq\frac{ C}{n^{d+2} (\log n)^2}\cdot n^d \cdot
n^2\rightarrow
0,
\]
and since $\operatorname{Var}_\mathbb{P}(\widetilde\mu _x^2) \le c n^6$,
\[
\operatorname{Var}_\mathbb{P}(B_n)\leq\frac{ C}{n^{2d+4} (\log
n)^4} \cdot n^d\cdot
n^6\leq\frac{C}{n^{d-2}(\log n)^4}.
\]
Since this bound is summable, \eqref{e:claim} follows.

It remains to show that for any $\delta_1>0$, $\mathbb{P}$-a.s.,
\begin{eqnarray*}
&&\limsup_n \frac{2}{(\log n)^2} \sum_{|x|,|y|\leq Kn}
\widetilde\mu _x\widetilde\mu _y\cdot\int_\delta^{t} p^\omega_{n^2 s}(0,x)
\int_{\delta_1}^{t-s}p^\omega_{n^2 r}(x,y) \, dr\,   ds\\
&&\qquad\leq 8(1+\varepsilon) \int_{|x|,|y|\leq K}  \biggl(\int_\delta^{t}
k_{ s}(0,x)\int_{0}^{t-s}k_{r}(x,y) \, dr\, ds \biggr)\,  dx\, dy.
\end{eqnarray*}
This follows easily from Theorem~\ref{thm:unif_llt} and Lemma~\ref{lem:DCT_bd}.
\end{pf}

\begin{pf*}{Proof of Theorem~\ref{thm:conv}}
By Lemma~\ref{lem:skorohod_marginal}, it suffices to show that for
any $t>0$ and $0<\varepsilon <t/2$, for $\mathbb{P}$-a.a.
${\omega}$,
%
\begin{equation}\label{eq:conv_prob}
\lim_n P^0_{\omega}\bigl(\big| S^{(n)}_t - 2t \big| \geq\epsilon\bigr)\leq\epsilon.
\end{equation}
Write
%
\begin{eqnarray}\label{e:Sndec}
\hspace*{35pt}
S^{(n)}_t - 2t &=& \bigl(S^{(n)}_t - \widetilde S^{(n)}_t \bigr)
+ \widetilde S^{(n)}_\delta
+ \bigl(\widetilde S^{(n)}_t - \widetilde S^{(n)}_\delta
- E^0_{\omega} \bigl(\widetilde S^{(n)}_t - \widetilde S^{(n)}_\delta \bigr) \bigr)\nonumber
\\[-8pt]\\[-8pt]
&&{} +  \bigl( E^0_{\omega} \bigl(\widetilde S^{(n)}_t - \widetilde S^{(n)}_\delta \bigr)
- 2 A_1(K, t, \delta) \bigr) + \bigl( 2 A_1(K, t, \delta) - 2t\bigr).\nonumber
\end{eqnarray}
By Proposition~\ref{prop:conv_exp}, $\mathbb{P}$-a.s., $ (
E^0_{\omega} (\widetilde S^{(n)}_t - \widetilde S^{(n)}_\delta
 )- 2 A_1(K, t,
\delta) )\rightarrow 0$. Let $0<\varepsilon _0<\varepsilon
/16$, to be chosen later.
Choose $K$ large enough so that the LHS in \eqref{ineq:big_ball}
is bounded by $\varepsilon _0$, and also
%
\begin{equation}\label{e:exp_diff}
\sup_{0< \delta\le t} |A_1(K,t,\delta) - (t-\delta) | \le
\varepsilon _0<\varepsilon /16.
\end{equation}
Now choose $a>0$ large enough so that the LHS in
\eqref{ineq:big_weight} is also bounded by~$\varepsilon _0$. Hence, for
all large $n$,
\[
P^0_{\omega}\bigl( \big|S^{(n)}_t - \widetilde S^{(n)}_t\big| >0 \bigr) \le2
\varepsilon _0< \varepsilon /4.
\]
Next choose $0<\delta<t/2$ so that by Lemma~\ref{lem:ini_ignore}
for all sufficiently large~$n$, $E^0_{\omega}\widetilde S^{(n)}_\delta
< \varepsilon ^2/16$,
and hence $P^0_{\omega}( \widetilde S^{(n)}_\delta
> \varepsilon /4) \le\varepsilon /4$.
Furthermore, by Propositions~\ref{prop:conv_exp} and
\ref{prop:conv_var_qu} and \eqref{e:exp_diff},
\begin{eqnarray*}
\limsup_n \operatorname{Var}_\mathbb{P}
\bigl(\widetilde S^{(n)}_t -\widetilde S^{(n)}_\delta \bigr)
&\leq&\varepsilon _0 + 8(1+\varepsilon _0)\cdot(t-\delta)^2/2 -
\bigl(2(t-\delta-\varepsilon _0)\bigr)^2\\
&\leq&\varepsilon _0 ( 1 + 4t^2 + 4t);
\end{eqnarray*}
hence by Chebyshev's inequality,
\[
\limsup_n P^0_{\omega}\bigl( \big|\widetilde S^{(n)}_t - \widetilde
S^{(n)}_\delta- E^0_{\omega}\bigl(\widetilde S^{(n)}_t - \widetilde S^{(n)}
_\delta\bigr)\big| \ge\varepsilon /4\bigr)
\leq16( 1 + 4t^2 + 4t)\cdot\varepsilon _0/\varepsilon ^2.
\]
Taking $\varepsilon _0$ so small that $\varepsilon _0<\varepsilon
/16$ and $16( 1 + 4t^2 +
4t)\cdot\varepsilon _0/\varepsilon ^2\leq\varepsilon /4$, we
obtain \eqref{eq:conv_prob}.
\end{pf*}

%

\printaddresses


\begin{thebibliography}{11}

\bibitem{BD08}
%
\begin{barticle}[vtex]
\bauthor{\bsnm{Barlow},~\bfnm{M.~T.}\binits{M.~T.}} \AND
\bauthor{\bsnm{Deuschel},~\bfnm{J.-D.}\binits{J.-D.}}
(\byear{2010}).
\btitle{Invariance principle for the random conductance model with
unbounded conductances}.
\bjournal{Ann. Probab.}
\bvolume{38}
\bpages{234--276}.
\end{barticle}
%
\endbibitem

\bibitem{BH09}
%
\begin{barticle}[vtex]
\bauthor{\bsnm{Barlow},~\bfnm{M.~T.}\binits{M.~T.}} \AND
\bauthor{\bsnm{Hambly},~\bfnm{B.~M.}\binits{B.~M.}}
(\byear{2009}).
\btitle{Parabolic {H}arnack inequality and local limit theorem for percolation
clusters}.
\bjournal{Electron. J. Probab.}
\bvolume{14}
\bpages{1--27}.
\bid{mr={2471657}}
\end{barticle}
%
\endbibitem

\bibitem{BC09}
%
\begin{bmisc}[vtex]
\bauthor{\bsnm{Barlow},~\bfnm{M.~T.}\binits{M.~T.}} \AND
\bauthor{\bsnm{{\v{C}}ern{\'y}},~\bfnm{Ji{\v{r}}{\'{\i}}}\binits{J.}}
(\byear{2009}).
Convergence to fractional kinetics for random walks associated with
unbounded conductances.
Preprint.
\end{bmisc}
%
\endbibitem

\bibitem{BC07}
%
\begin{barticle}[mr]
\bauthor{\bsnm{Ben~Arous},~\bfnm{G{\'e}rard}\binits{G.}} \AND
\bauthor{\bsnm{{\v{C}}ern{\'y}},~\bfnm{Ji{\v{r}}{\'{\i}}}\binits{J.}}
(\byear{2007}).
\btitle{Scaling limit for trap models on {$\Bbb Z\sp d$}}.
\bjournal{Ann. Probab.}
\bvolume{35}
\bpages{2356--2384}.
\bid{mr={2353391}}
\end{barticle}
%
\endbibitem

\bibitem{BC08}
%
\begin{barticle}[mr]
\bauthor{\bsnm{Ben~Arous},~\bfnm{G{\'e}rard}\binits{G.}} \AND
\bauthor{\bsnm{{\v{C}}ern{\'y}},~\bfnm{Ji{\v{r}}{\'{\i}}}\binits{J.}}
(\byear{2008}).
\btitle{The arcsine law as a universal aging scheme for trap models}.
\bjournal{Comm. Pure Appl. Math.}
\bvolume{61}
\bpages{289--329}.
\bid{mr={2376843}}
\end{barticle}
%
\endbibitem

\bibitem{BCM06}
%
\begin{barticle}[mr]
\bauthor{\bsnm{Ben~Arous},~\bfnm{G{\'e}rard}\binits{G.}},
\bauthor{\bsnm{{\v{C}}ern{\'y}},~\bfnm{Ji{\v{r}}{\'{\i}}}\binits
{J.}} \AND
\bauthor{\bsnm{Mountford},~\bfnm{Thomas}\binits{T.}}
(\byear{2006}).
\btitle{Aging in two-dimensional {B}ouchaud's model}.
\bjournal{Probab. Theory Related Fields}
\bvolume{134}
\bpages{1--43}.
\bid{mr={2221784}}
\end{barticle}
%
\endbibitem

\bibitem{BP07}
%
\begin{barticle}[vtex]
\bauthor{\bsnm{Biskup},~\bfnm{Marek}\binits{M.}} \AND
\bauthor{\bsnm{Prescott},~\bfnm{Timothy~M.}\binits{T.~M.}}
(\byear{2007}).
\btitle{Functional {CLT} for random walk among bounded random conductances}.
\bjournal{Electron. J. Probab.}
\bvolume{12}
\bpages{1323--1348 (electronic)}.
\bid{mr={2354160}}
\end{barticle}
%
\endbibitem

\bibitem{Grim}
%
\begin{bbook}[mr]
\bauthor{\bsnm{Grimmett},~\bfnm{Geoffrey}\binits{G.}}
(\byear{1999}).
\btitle{Percolation},
\bedition{2nd} ed.
\bseries{Grundlehren der Mathematischen Wissenschaften [Fundamental Principles
of Mathematical Sciences]}
\bvolume{321}.
\bpublisher{Springer}, \baddress{Berlin}.
\bid{mr={1707339}}
\end{bbook}
%
\endbibitem

\bibitem{JS03}
%
\begin{bbook}[mr]
\bauthor{\bsnm{Jacod},~\bfnm{Jean}\binits{J.}} \AND
\bauthor{\bsnm{Shiryaev},~\bfnm{Albert~N.}\binits{A.~N.}}
(\byear{2003}).
\btitle{Limit Theorems for Stochastic Processes},
\bedition{2nd} ed.
\bseries{Grundlehren der Mathematischen Wissenschaften [Fundamental Principles
of Mathematical Sciences]}
\bvolume{288}.
\bpublisher{Springer}, \baddress{Berlin}.
\bid{mr={1943877}}
\end{bbook}
%
\endbibitem

\bibitem{LPW}
%
\begin{bbook}[vtex]
\bauthor{\bsnm{Levin},~\bfnm{David~A.}\binits{D.~A.}},
\bauthor{\bsnm{Peres},~\bfnm{Yuval}\binits{Y.}} \AND
\bauthor{\bsnm{Wilmer},~\bfnm{Elizabeth~L.}\binits{E.~L.}}
(\byear{2009}).
\btitle{Markov Chains and Mixing Times}.
\bpublisher{Amer. Math. Soc.}, \baddress{Providence, RI}.
\bid{mr={2466937}}
\end{bbook}
%
\endbibitem

\bibitem{Mat08}
%
\begin{barticle}[mr]
\bauthor{\bsnm{Mathieu},~\bfnm{P.}\binits{P.}}
(\byear{2008}).
\btitle{Quenched invariance principles for random walks with random
conductances}.
\bjournal{J. Stat. Phys.}
\bvolume{130}
\bpages{1025--1046}.
\bid{mr={2384074}}
\end{barticle}
%
\endbibitem

\end{thebibliography}
\end{document}